\documentclass[12pt]{article}
\usepackage{amsmath,amsthm,amsfonts,amssymb}
\usepackage{fullpage}
\usepackage{multirow}
\usepackage{array}
\usepackage{bbm}
\usepackage[noblocks]{authblk}
\usepackage{verbatim}
\usepackage{tikz}
\usepackage{float}
\usepackage{enumerate}
\usepackage{diagbox}
\usepackage{soul}
\usepackage{url}
\usepackage{mathtools}
\usepackage{hyperref}
\usepackage{listings}

\newtheorem{theorem}{Theorem}[section]

\newtheorem{corollary}[theorem]{Corollary}

\theoremstyle{definition}

\newlength{\Oldarrayrulewidth}

\begin{document}

\author{Hailey~Evans\thanks{hmevans@cedarcrest.edu}}
\author{Joshua~Harrington\thanks{joshua.harrington@cedarcrest.edu}}
\author{Kendall~Heiney\thanks{KMHeiney@cedarcrest.edu}}
\author{Maggie~Wieczorek\thanks{margaretwieczorek@cedarcrest.edu}}
\affil{Department of Mathematics, Cedar Crest College}
\date{\today}

\title{Generalized Sierpi\'{n}ski and Riesel numbers of the form $tb^t+\alpha$}
\maketitle

\begin{abstract}
Let $b\geq 2$ be an integer. We call an integer $k$ a $b$-Sierpi\'{n}ski number if $\gcd(k+1,b-1)=1$ and $k\cdot b^n+1$ is composite for all positive integers $n$. We similarly call $k$ a $b$-Riesel number if $\gcd(k-1,b-1)=1$ and $k\cdot b^n-1$ is composite for all positive integers $n$. An integer that is simultaneously $b$-Sierpi\'{n}ski and $b$-Riesel is called a $b$-Brier number. In this article, we show that for any integer $\alpha\neq 0$, there are infinitely many $b$-Sierpi\'{n}ski numbers and infinitely many $b$-Riesel numbers of the form $tb^t+\alpha$. We further show that when $b+1$ is not a power of $2$, there are infinitely $b$-Brier number of this form.

\end{abstract}

\section{Introduction}
In 1956, Riesel showed that there are infinitely many odd integers $k$ such that $k\cdot 2^n-1$ is composite for all positive integers $n$ \cite{riesel}. Today, such numbers are known as Riesel numbers. Four years after Riesel's paper, Sierpi\'{n}ski showed that there are infinitely many odd integers $k$ such that $k\cdot 2^n+1$ is composite for all positive integers $n$ \cite{sierpinski}; these numbers are now known as Sierpi\'{n}ski numbers. A number that is both Riesel and Sierpi\'{n}ski is known as a Brier number.

The definitions of Riesel and Sierpi\'{n}ski numbers can be generalized in the following way: for an integer $b\geq 2$, a positive integer $k$ with $\gcd(k-1,b-1)=1$ is called a $b$-Riesel number if $k\cdot b^n-1$ is composite for all positive integers $n$. Similarly, $k$ is called a $b$-Sierpi\'{n}ski number if $\gcd(k+1,b-1)=1$ and $k\cdot b^n+1$ is composite for all positive integers $n$. We analogously define a $b$-Brier number to be a number that is both $b$-Riesel and $b$-Sierpi\'{n}ski. Note that the requirement that $\gcd(k-\beta,b-1)$, with $\beta\in\{-1,1\}$, is to avoid the trivial case where every term of $k\cdot b^n+\beta$ has a common divisor greater than $1$. 

Since the original findings of Riesel and Sierpi\'{n}ski, many papers have been published that investigate the existence of Riesel numbers,  Sierpi\'{n}ski numbers, and Brier numbers in various other integer sequences, including binomial coefficients \cite{abbdehsw}, polygonal numbers \cite{befks,be}, the Lucas sequence \cite{bff}, Carmichael numbers \cite{bflps}, perfect powers \cite{chen, ffk, fj}, Ruth-Aaron pairs \cite{efk}, the image of various polynomials \cite{fhj}, the Fibonacci sequence \cite{lm}, and Narayana's cow sequence \cite{fh}. Of note to this article, in 2012, Berrizbeitia et al. studied the existence of Brier numbers in the sequence of Cullen numbers and the sequence of Woodall numbers \cite{bfglm}.

A Cullen number is one of the form $C(t)=t\cdot 2^t+1$, where $t$ is a non-negative integer. These numbers were first studied by James Cullen in 1905. Similarly, numbers of the form $W(t)=t\cdot 2^t-1$, introduced by  Allan J. C. Cunningham and H. J. Woodall, are called Woodall numbers. For an integer $b\geq 2$, we can generalize each of these sequences by defining $C_b(t)=t\cdot b^t+1$ and $W_b(t)=t\cdot b^t-1$ to be a $b$-Cullen number and a $b$-Woodall number, respectively. These sequences have traditionally been studied in the pursuit of prime numbers. 

Berrizbeitia et al. \cite{bfglm}, showed that the Cullen number $C(t)$ is a Brier number for any $t$ satisfying
$$t\equiv 41 988 231 525 567 017 419 567 132 158 576\pmod{51 409 261 600 541 515 243 210 608 930 960},$$
and the Woodall number $W(t)$ is a Brier number for any $t$ satisfying
$$t\equiv 10 922 075 158 847 048 672 240 898 358 992
\pmod{51 409 261 600 541 515 243 210 608 930 960}.$$

In this article, we extend these results in several ways. First, we show that for any integers $b\geq 2$ and $\alpha$ and $\beta$ with $\alpha\beta\neq 0$, there exist infinitely many integers $t$ with $\gcd((tb^t+\alpha)+\beta,b-1)=1$ such that $(tb^t+\alpha)b^n+\beta$ is composite for all positive integers $n$. We then show that when $b+1$ is not a power of $2$, there are infinitely many integers $t$ for which $tb^t+\alpha$ is a $b$-Brier number.

\section{The Covering System Method}
The common approach used to find Riesel numbers and Sierpi\'{n}ski numbers involves the use of covering systems of the integers, a concept that was first introduced by Erd\H{o}s \cite{erdos}. A covering system of the integers, which we often shorten to covering system, is a finite collection of congruences such that every integer satisfies at least one congruence in the collection. For example, the congruences $0\pmod{2}$ and $1\pmod{2}$ form a covering system since every integer is either even or odd. Recognizing that every odd integer is either $1\pmod{4}$ or $3\pmod{4}$, the collection 
\begin{align*}
0\pmod{2}\\
1\pmod{4}\\
3\pmod{4}
\end{align*}
forms a slightly less trivial example of a covering system.

In practice, it is often useful to create covering systems whose moduli are distinct. Following the reasoning above, to construct a covering system with distinct moduli, we can start with the collection of congruences
\begin{eqnarray}\label{eq:powersof2moduli}
0\pmod{2}\nonumber\\
1\pmod{4}\nonumber\\
3\pmod{8}\\
\vdots\qquad\qquad\nonumber\\
2^{m-1}-1\pmod{2^m},\nonumber
\end{eqnarray}
where $m\geq 2$ is an integer.
Notice that every integer $n\not\equiv 2^m-1\pmod{2^m}$ will satisfy exactly one of the congruences in \eqref{eq:powersof2moduli}. Let $q\leq m+1$ be an odd prime and, using the Chinese remainder theorem, let $a_j$ satisfy the two congruences
\begin{align*}
a_j&\equiv j\pmod{q}\\
a_j&\equiv 2^m-1\pmod{2^j}
\end{align*} 
for $0\leq j\leq q-1$. Since every integer satisfies one of the congruences in the set $\{n\equiv j\pmod{q}:0\leq j\leq q-1\}$, every integer satisfying $n\equiv 2^m-1\pmod{2^m}$ must satisfy one of the congruences in the set $\{n\equiv a_j\pmod{q\cdot 2^j}:0\leq j\leq q-1\}$. Hence, for any integer $m\geq 2$ and odd prime $q\leq m+1$, the collection of congruences 
\begin{equation}\label{eq:firstcovering}
\begin{aligned}
2^{j-1}-1\pmod{2^j}&\text{ for }1\leq j\leq q-1\\
a_j\pmod{q\cdot 2^j}&\text{ for }0\leq j\leq q-1
\end{aligned}
\end{equation}
forms a covering system. 

To construct Sierpi\'{n}ski and Riesel numbers, we will make use of the following theorem, which is a special case of a theorem of Zsigmondy \cite{zsigmondy}.
\begin{theorem}\label{thm:zsigmondy}
Let $b\geq 2$ be an integer and let $w$ be a positive integer. Then there exists a prime $p$ that divides $b^w-1$ and does not divide $b^{\ell}-1$ for any $\ell<w$ except in the following cases:
\begin{itemize}
\item $b=2$ and $w\in\{1,6\}$
\item $b+1$ is a power of $2$ and $w=2$.
\end{itemize} 
\end{theorem}
For positive integers $w$ and $b\geq 2$, we call a prime $p$ a primitive prime divisor of $b^w-1$ if $p$ divides $b^w-1$ and $p$ does not divide $b^{\ell}-1$ for all $1\leq\ell<w$. 

Let $\mathcal{C}=\{r_i\pmod{m_i}:1\leq i\leq \nu\}$ be a covering system of the form given in \eqref{eq:firstcovering} with $q\geq 5$. We may assume that $m_\iota<m_i$ for any $\iota<i$. Note that since $q\geq 5$, $m=6$ is not a modulus in $\mathcal{C}$. Thus, when $b+1$ is not a power of $2$, for each $1\leq i\leq \nu$ we can use Theorem~\ref{thm:zsigmondy} to let $p_i$ be a primitive prime divisor of $b^{m_i}-1$. We note that there may be multiple primitive prime divisors of $b^{m_i}-1$, and we may choose any one of them for our purposes. Choosing $p_i$ in this way ensures that $p_1,\ldots,p_{\nu}$ are distinct. Then for any integer $k$, when $n$ is positive and $n\equiv r_i\pmod{m_i}$, we have that
$$k\cdot b^n+1\equiv k\cdot b^{r_i}+1\pmod{p_i}.$$
Consequently, letting $k>p_i$ and $k\equiv -b^{-r_i}\pmod{p_i}$ ensures that $k\cdot b^n+1$ is divisible by $p_i$, and, hence, is composite when $n\equiv r_i\pmod{m_i}$. Since each $p_i$ is distinct, we can use the Chinese remainder theorem to let $k$ satisfy the congruences $k\equiv -b^{-r_i}\pmod{p_i}$ for each $1\leq i\leq \nu$. Since $\mathcal{C}$ is a covering system, this choice of $k$ guarantees that $k\cdot b^n+1$ is composite for all positive integers $n$ and is, therefore, a $b$-Sierpi\'{n}ski number.

In light of the process described above, we say that a covering system $\mathcal{C}=\{r_i\pmod{m_i}:1\leq i\leq \nu\}$ is a $b$-primitive covering system if there exist distinct primes $p_1,p_2,\ldots,p_\nu$ such that for each $1\leq i\leq \nu$, $p_i$ is a primitive prime divisor of $b^{m_i}-1$. Notice that the covering system provided in \eqref{eq:firstcovering} is a $b$-primitive covering when $q\geq 5$ and $b+1$ is not a power of $2$. Throughout the paper, for a given $b$-primitive covering system $\mathcal{C}$, we let $\mathcal{P}_b(\mathcal{C})=\{p_1,p_2,\ldots,p_\nu\}$ be a set of primitive prime divisors chosen as we have described here. 

\section{Main Results}
The first theorem of this section utilizes the covering system given in \eqref{eq:firstcovering} with $m=q-1$. As a consequence of Theorem~\ref{thm:zsigmondy}, we have the restriction $b+1$ is not a power of $2$ in this theorem.
\begin{theorem}\label{thm:notbrierwith2}
Let $b\geq 2$ be an integer with $b+1$ not a power of $2$ and let $\alpha$ and $\beta$ be integers with $\alpha\beta\neq 0$. Further, let $q\geq 5$ be a prime and let $\mathcal{C}$ be a $(b,1)$-primitive covering system with moduli $m_{(j,0)}=2^j$ for $1\leq j\leq q-1$ and $m_{(j,1)}=2^jq$ for $0\leq j\leq q-1$. Lastly, let $M=\text{lcm}\left(\{2^{q-1},b-1,q\}\cup\mathcal{P}_b(\mathcal{C})\right)$. There exists an integer $T$ such that for all $t\equiv T\pmod{M}$ with $(tb^t+\alpha)b+\beta\geq\max\{\mathcal{P}_b(\mathcal{C})\}$, $\gcd(tb^t+\alpha+\beta,b-1)=1$, and the expression $(tb^t+\alpha)b^n+\beta$ is composite for all positive integers $n$.
\end{theorem}
\begin{proof}
Let $\mathcal{C}=\{r_{(j,0)}\pmod{m_{(j,0)}}:1\leq j\leq q-1\}\cup\{r_{(j,1)}\pmod{m_{(j,1)}}:0\leq j\leq q-1\}$ and let $p_{(j,j')}$ be a primitive prime divisor of $b^{m_{(j,j')}}-1$ when $0\leq j\leq q-1$ and $0\leq j'\leq 1$ with $(j,j')\neq(0,0)$. Note that since the moduli of $\mathcal{C}$ are distinct, the primes in $\mathcal{P}_b(\mathcal{C})$ are distinct and relatively prime to $b-1$. To simplify our arguments, we define $P$ to be the product of the primes in $\mathcal{P}(\mathcal{C})$, and we let $r_{(0,0)}=m_{(0,0)}=p_{(0,0)}=1$.  

Notice that if $2$ divides $b^m-1$, then $2$ divides $b-1$, and hence, since all moduli of $\mathcal{C}$ are greater than $1$, $2\notin\mathcal{P}_b(\mathcal{C})$. Since $b^{2^{q-1}q}\equiv b^{2^{q-1}}\pmod{q}$, we have that if $q$ divides $P$, then $q=p_{(a,0)}$ for some $1\leq a\leq q-1$. For the remainder of the proof, if $q$ divides $P$, then we define $a$ as mentioned here. Since all moduli of $\mathcal{C}$ are greater then $1$, if $q$ divides $P$, then $q$ does not divide $b-1$, and, therefore, $\gcd(q,b-1)=1$.

With the above observations, we can use the Chinese remainder theorem to let $\tau$ satisfy the congruences 
\begin{align*}
        \tau&\equiv 1-\alpha-\beta\pmod{2^{q-1}}\\
        \tau&\equiv\begin{cases} 1-\alpha-\beta\pmod{q}&\text{ if }q\nmid P\\(-\alpha-\beta\cdot b^{-r_{(a,0)}})\cdot b^{-1+\alpha+\beta}\pmod{q}&\text{ if }q\mid P,\end{cases}
    \end{align*}
let $\mathcal{T}$ satisfy the congruences    \begin{equation*}
    \begin{aligned}
        \mathcal{T}&\equiv1-\alpha-\beta\pmod{2^{q-1}}\\
        \mathcal{T}&\equiv1-\alpha-\beta\pmod{b-1}\\
        \mathcal{T}&\equiv(-\alpha-\beta\cdot b^{-r_{(j,j')}})\cdot b^{-\tau}\pmod{p_{(j,j')}}&\text{ for }0\leq j\leq q\text{ and }0\leq j'\leq 1 \\
        &&\text{ with }p_{(j,j')}\neq q,
    \end{aligned}
    \end{equation*}
and let $T$ satisfy the congruences 
    \begin{align*}
        T&\equiv\mathcal{T}\pmod{(b-1)\cdot P/\gcd(q,P)}\\
        T&\equiv\tau\pmod{2^{q-1}q}.
    \end{align*}
Then we can write $T=2^{q-1}qh+\tau$ for some integer $h$, and since $b^{2^{q-1}qh}\equiv 1\pmod{p_{(j,j')}}$ for each $0\leq j\leq q-1$ and $0\leq j'\leq 1$, we have that when $p_{(j,j')}\neq q$, 
    \begin{align*}
        T\cdot b^T+\alpha
        &=T\cdot b^{2^{q-1}qh}\cdot b^{\tau}+\alpha\\
        &\equiv ((-\alpha-\beta\cdot b^{-r_{(j,j')}})\cdot b^{-\tau})\cdot b^{2^{q-1}qh}\cdot b^{\tau}+\alpha\pmod{p_{(j,j')}}\\
        &\equiv-\beta\cdot b^{-r_{(j,j')}}\pmod{p_{(j,j')}}.
    \end{align*}

Suppose that $q$ divides $P$. Then $q=p_{(a,0)}$ and $T\equiv\tau\pmod{2^{q-1}q}$ imply $T\equiv-1-\alpha-\beta\pmod{2^{q-1}}$. Thus, $T=2^{q-1}\ell+1-\alpha-\beta$ for some integer $\ell$, and since $b^{2^{q-1}\ell}\equiv 1\pmod{p_{(a,0)}}$, we have that 
    \begin{align*}
        T\cdot b^T+\alpha
        &=T\cdot b^{2^{q-1}\ell}\cdot b^{1-\alpha-\beta}+\alpha\\
        &\equiv((-\alpha-\beta\cdot b^{-r_{(a,0)}})\cdot b^{-1+\alpha+\beta})\cdot b^{2^{q-1}\ell}\cdot b^{1-\alpha-\beta}+\alpha\pmod{q}\\
        &\equiv-\beta\cdot b^{-r_{(a,0)}}\pmod{q}.
    \end{align*}
In all cases, if $n\equiv r_{(j,j')}\pmod{m_{j,j'}}$ for some $r_{(j,j')}\pmod{m_{(j,j')}}$ in $\mathcal{C}$, then $(Tb^T+\alpha)\cdot b^n+\beta$ is divisible by $p_{(j,j')}$. The condition that $T\equiv 1-\alpha-\beta\pmod{b-1}$ ensures that $(Tb^T+\alpha)+\beta\equiv 1\pmod{b-1}$, and, therefore, $\gcd((Tb^T+\alpha)+\beta,b-1)=1$. Since $\mathcal{C}$ is a covering system, this finishes the proof.
\end{proof}

When $\alpha\in\{-1,1\}$ and $\beta=\{-1,1\}$, Theorem~\ref{thm:notbrierwith2} establishes infinitely many Riesel numbers and Sierpi\'{n}ski numbers in the Cullen sequence and the Woodall sequence that were not previously given by Berrizbeitia et al. We present the next corollary to demonstrate the construction given in the proof of Theorem~\ref{thm:notbrierwith2}.
\begin{corollary}
Each of the following hold with $M=3868562622766813359059760$.
\begin{enumerate}
    \item\label{item:corcullensierpinski} For all $t\equiv2245377406103792702454767\pmod{M}$, the Cullen number $C(t)$ is a Sierpi\'{n}ski number.
    \item\label{item:corcullenriesel}For all $t\equiv
2215074033447763254589281\pmod{M}$, the Cullen number $C(t)$ is a Riesel number.
    \item\label{item:corwoodallsierpinski}For all $t\equiv
1951609044082776021493089\pmod{M}$, the Woodall number $W(t)$ is a Sierpi\'{n}ski number.
    \item\label{item:corwoodallriesel}For all $t\equiv
3334297893475587915471523\pmod{M}$, the Woodall number $W(t)$ is a Riesel number.
\end{enumerate}
\end{corollary}
\begin{proof}
Let $b=2$, $q=5$, and let 
    \begin{align*}
        \mathcal{C}=\{&0\pmod{2}, 1\pmod{4},3\pmod{8}, 7\pmod{16}\\
        &0\pmod{5}, 1\pmod{10}, 7\pmod{20}, 23\pmod{40},79\pmod{80}\}.
    \end{align*}
The corollary follows from Theorem~\ref{thm:notbrierwith2} by appropriately calculating $\tau$, $\mathcal{T}$, and $T$ in each of the four cases $(\alpha,\beta)\in\{(1,1),(1,-1),(-1,1),(-1,-1)\}$. We only present the construction for part~\ref{item:corcullensierpinski} of the corollary, the case when $(\alpha,\beta)=(1,1)$.

We find $\mathcal{P}_b(\mathcal{C})=\{3,5,11,17,31,257,41, 61681,4278255361\}$ and compute $$P=241785163922925834941235.$$ 
Then
    \begin{align*}
        \left\{\left(r_{(j,j')},m_{(j,j')},p_{(j,j')}\right):0\leq j\leq q-1, 0\leq j'\leq 1\right\}=&
        \{(0,2,3),(1,4,5),(3,8,17),(7,16,257),\\
        &(0,5,31),(1,10,11),(7,20,41),\\
        &(23,40,61681),(79,80,4278255361)\}.
    \end{align*}
Since $q$ divides $P$, we let $a=2$ and  $\left(r_{(a,0)},m_{(a,0)},p_{(a,0)}\right)=(1,4,5)$.

Following the proof of Theorem~\ref{thm:notbrierwith2}, to establish part~\ref{item:corcullensierpinski} of the corollary, we let $\tau=47$ so that $\tau$ satisfies the congruences
    \begin{align*}
            \tau&\equiv -1\pmod{16}\\
            \tau&\equiv (-1-2^{-1})\cdot 2\pmod{5}\equiv 2\pmod{5}.
    \end{align*}
We further let $\mathcal{T}=
697952356997067358830863$ so that $\mathcal{T}$ satisfies the congruences
    \begin{align*}
        \mathcal{T}&\equiv -1\pmod{16}\\
        \mathcal{T}&\equiv (-1-2^{-0})\cdot 2^{-47}\pmod{3}\equiv2\pmod{3}\\
        \mathcal{T}&\equiv (-1-2^{-3})\cdot 2^{-47}\pmod{17}\equiv2\pmod{17}\\
        \mathcal{T}&\equiv (-1-2^{-7})\cdot 2^{-47}\pmod{257}\equiv2\pmod{257}\\
        \mathcal{T}&\equiv (-1-2^{-0})\cdot 2^{-47}\pmod{31}\equiv15\pmod{31}\\
        \mathcal{T}&\equiv (-1-2^{-1})\cdot 2^{-47}\pmod{11}\equiv10\pmod{11}\\
        \mathcal{T}&\equiv (-1-2^{-7})\cdot 2^{-47}\pmod{41}\equiv26\pmod{41}\\
        \mathcal{T}&\equiv (-1-2^{-23})\cdot 2^{-47}\pmod{61681}\equiv7168\pmod{61681}\\
        \mathcal{T}&\equiv (-1-2^{-79})\cdot 2^{-47}\pmod{4278255361}\equiv4177983751\pmod{4278255361}.      
    \end{align*}
Finally, we let $T=2245377406103792702454767$ so that $T$ satisfies the congruences
    \begin{align*}
        T&\equiv 20953898012875020995405\pmod{ 48357032784585166988247}\\
        T&\equiv47\pmod{80}.
    \end{align*}
\end{proof}

Our next theorem removes the restriction that $b+1$ is not a power of $2$ by using a covering system that does not have $2$ as a modulus. An example of such a covering system is presented in \eqref{eq:lastcovering} at the end of this section. Since the covering system in \eqref{eq:lastcovering} has $6$ as a modulus, as a consequence of Theorem~\ref{thm:zsigmondy}, we have the restriction $b\geq 3$ in our next theorem.

\begin{theorem}\label{thm:notbrierwithout2}
Let $b\geq 3$ be an integer and let $\alpha\beta\neq 0$. Further, let $q\geq 5$ be a prime and let $\mathcal{C}$ be a $(b,1)$-primitive covering system with moduli $m_{(j,\overline{j},j')}=2^j3^{\overline{j}}q^{j'}$ for $0\leq j\leq q$, $1\leq\overline{j}\leq q$, and $0\leq j'\leq 1$. Lastly, let $M=\text{lcm}\left(\{2^q,3^q,b-1,q\}\cup\mathcal{P}_b(\mathcal{C})\right)$. There exists an integer $T$ such that for all $t\equiv T\pmod{M}$ with $(tb^t+\alpha)b+\beta\geq\max\{\mathcal{P}_b(\mathcal{C})\}$, $\gcd(tb^t+\alpha+\beta,b-1)=1$, and the expression $(tb^t+\alpha)b^n+\beta$ is composite for all positive integers $n$.
\end{theorem}
\begin{proof}
Let $\mathcal{C}=\{r_{(j,\overline{j},j')}\pmod{m_{(j,\overline{j},j')}}:0\leq j\leq q,1\leq\overline{j}\leq q, 0\leq j'\leq 1\}$ and let $p_{(j,\overline{j},j')}$ be a primitive prime divisor of $b^{m_{(j,\overline{j},j')}}-1$ when $0\leq j\leq q$, $1\leq\overline{j}\leq q$, and $0\leq j'\leq 1$. Note that since the moduli of $\mathcal{C}$ are distinct, the primes in $\mathcal{P}_b(\mathcal{C})$ are distinct and relatively prime to $b-1$. To simplify our arguments, we define $P$ to be the product of the primes in $\mathcal{P}(\mathcal{C})$, and we let $r_{(0,0,0)}=m_{(0,0,0)}=p_{(0,0,0)}=1$.  

Notice that if $b$ is divisible by $3$, then $3\notin\mathcal{P}_b(\mathcal{C})$. If $b$ is not divisible by $3$, then $b^2-1$ is divisible by $3$ by Fermat's little theorem. Consequently, since all moduli of $\mathcal{C}$ are greater than 2, $3\notin\mathcal{P}_b(\mathcal{C})$. Since $b^{2^q3^qq}\equiv b^{2^q3^q}\pmod{q}$, we have that if $q$ divides $P$, then $q=p_{(a,\overline{a},0)}$ for some $0\leq a\leq q$ and $1\leq\overline{a}\leq q$. For the remainder of the proof, if $q$ divides $P$, then we define $a$ and $\overline{a}$ as mentioned here. Since all moduli of $\mathcal{C}$ are greater then $1$, if $q$ divides $P$, then $q$ does not divide $b-1$, and, therefore, $\gcd(q,b-1)=1$.

With the above observations, we can use the Chinese remainder theorem to let $\tau$ satisfy the congruences
    \begin{align*}
        \tau&\equiv 1-\alpha-\beta\pmod{2^q3^q}\\
        \tau&\equiv\begin{cases} 1-\alpha-\beta\pmod{q}&\text{ if }q\nmid P\\(-\alpha-\beta\cdot b^{-r_{(a,\overline{a},0)}})\cdot b^{-1+\alpha+\beta}\pmod{q}&\text{ if }q\mid P,\end{cases}
    \end{align*}
let $\mathcal{T}$ satisfy the congruences 
    \begin{equation*}\label{eq:cullenwierpinskifirstcongruences}
    \begin{aligned}
        \mathcal{T}&\equiv1-\alpha-\beta\pmod{2^q3^q}\\
        \mathcal{T}&\equiv1-\alpha-\beta\pmod{b-1}\\
        \mathcal{T}&\equiv(-\alpha-\beta\cdot b^{-r_{(j,\overline{j},j')}})\cdot b^{-\tau}\pmod{p_{(j,\overline{j},j')}}&\text{ for }0\leq j\leq q\text{ and }0\leq j'\leq 1 \\
        &&\text{ with }p_{(j,\overline{j},j')}\neq q,
    \end{aligned}
    \end{equation*}
and let $T$ satisfy the congruences 
    \begin{align*}
        T&\equiv\mathcal{T}\pmod{(b-1)\cdot P/\gcd(q,P)}\\
        T&\equiv\tau\pmod{2^q3^qq}.
    \end{align*}
Then we can write $T=2^q3^qqh+\tau$ for some integer $h$, and since $b^{2^q3^qqh}\equiv 1\pmod{p_{(j,\overline{j},j')}}$ for each $0\leq j\leq q$, $1\leq\overline{j}\leq q$, and $0\leq j'\leq 1$, we have that when $p_{(j,\overline{j},j')}\neq q$, 
    \begin{align*}
        T\cdot b^T+\alpha
        &=T\cdot b^{2^q3^qqh}\cdot b^{\tau}+\alpha\\
        &\equiv ((-\alpha-\beta\cdot b^{-r_{(j,\overline{j},j')}})\cdot b^{-\tau})\cdot b^{2^q3^qqh}\cdot b^{\tau}+\alpha\pmod{p_{(j,\overline{j},j')}}\\
        &\equiv-\beta\cdot b^{-r_{(j,\overline{j},0)}}\pmod{p_{(j,\overline{j},j')}}.
    \end{align*}

Suppose that $q$ divides $P$. Then $q=p_{(a,\overline{a},0)}$ and $T\equiv\tau\pmod{2^q3^qq}$ imply $T\equiv-1-\alpha-\beta\pmod{2^q3^q}$. Thus, $T=2^q3^q\ell+1-\alpha-\beta$ for some integer $\ell$, and since $b^{2^q3^q\ell}\equiv 1\pmod{p_{(a,\overline{a},0)}}$, we have that 
    \begin{align*}
        T\cdot b^T+\alpha
        &=T\cdot b^{2^q3^q\ell}\cdot b^{1-\alpha-\beta}+\alpha\\
        &\equiv((-\alpha-\beta\cdot b^{-r_{(a,\overline{a},0)}})\cdot b^{-1+\alpha+\beta})\cdot b^{2^q3^q\ell}\cdot b^{1-\alpha-\beta}+\alpha\pmod{q}\\
        &\equiv-\beta\cdot b^{-r_{(a,\overline{a},0)}}\pmod{q}.
    \end{align*}
In all cases, if $n\equiv r_{(j,\overline{j},j')}\pmod{m_{j,\overline{j},j'}}$ for some $r_{(j,\overline{j},j')}\pmod{m_{(j,\overline{j},j')}}$ in $\mathcal{C}$, then $(T\cdot b^T+\alpha)\cdot b^n+\beta$ is divisible by $p_{(j,\overline{j},j')}$. The condition that $T\equiv 1-\alpha-\beta\pmod{b-1}$ ensures that $(Tb^T+\alpha)+\beta\equiv 1\pmod{b-1}$, and, therefore, $\gcd((Tb^T+\alpha)+\beta,b-1)=1$. Since $\mathcal{C}$ is a covering system, this finishes the proof.
\end{proof}

 By letting $\alpha\in\{-1,1\}$ and $\beta\in\{-1,1\}$, our next corollary follows immediately from Theorem~\ref{thm:notbrierwith2} and Theorem~\ref{thm:notbrierwithout2}, together with the covering systems provided in \eqref{eq:firstcovering} and \eqref{eq:lastcovering}.

\begin{corollary}
For any integer $b\geq 2$, there are infinitely many $b$-Sierpi\'{n}ski numbers and infinitely many $b$-Riesel numbers in the sequence of $b$-Cullen numbers and the sequence of $b$-Woodall numbers.
\end{corollary}

The final theorem of this paper utilizes the covering system given in \eqref{eq:firstcovering} and the covering system given in \eqref{eq:lastcovering}. As a consequence of Theorem~\ref{thm:zsigmondy}, we have the restrictions that $b\geq 3$ and $b$ not a power of $2$.

\begin{theorem}\label{thm:brier}
Let $b\geq 3$ be an integer with $b+1$ not a power of $2$ and let $\alpha$ be an integer with $\alpha\neq0$. Further, let $q\geq 5$ be a prime, let $\mathcal{C}_1$ be a $(b,1)$-primitive covering system with moduli $m_{(j,0,0)}=2^j$ for $1\leq j\leq q$ and $m_{(j,0,1)}=2^{j}q$ for $0\leq j\leq q$, and let $\mathcal{C}_2$ be a $(b,1)$-primitive covering system with moduli $m_{(j,\overline{j},j')}=2^j\cdot3^{\overline{j}}\cdot q^{j'}$ for each $0\leq j\leq q$, $1\leq \overline{j}\leq q$, and $0\leq j'\leq 1$. Lastly, let $M=\text{lcm}(\{2^q,3^q,b-1,q\}\cup\mathcal{P}_b(\mathcal{C}_1)\cup\mathcal{P}_b(\mathcal{C}_2)$. There exists an integer $T$ such that for all $t\equiv T\pmod{M}$ with $(tb^t+\alpha)b+1\geq\max\{\mathcal{P}_b(\mathcal{C}_1)\cup\mathcal{P}_b(\mathcal{C}_2\}$, $tb^t+\alpha$ is a $b$-Brier number.
\end{theorem}
\begin{proof}
Let $\mathcal{C}_1=\{r_{(j,0,0)}\pmod{m_{(j,0,0)}}:1\leq j\leq q\}\cup\{r_{(j,0,1)}\pmod{m_{(j,0,1)}}:0\leq j\leq q\}$ and let $p_{(j,0,j')}$ be a primitive prime divisor of $b^{m_{(j,0,j')}}-1$ when $0\leq j\leq q$ and $0\leq j'\leq 1$, with $(j,0,j')\neq(0,0,0)$. Similarly, let $\mathcal{C}_2=\{r_{(j,\overline{j},j')}\pmod{m_{(j,\overline{j},j')}}:0\leq j\leq q,1\leq \overline{j}\leq q,0\leq j'\leq 1\}$ and let $p_{(j,\overline{j},j')}$ be a primitive prime divisor of $b^{m_{(j,\overline{j},j')}}-1$ when $0\leq j\leq q$, $1\leq\overline{j}\leq q$, and $0\leq j'\leq 1$. Note that since the moduli of $\mathcal{C}_1$ and $\mathcal{C}_2$ are distinct, the primes in $\mathcal{P}_b(\mathcal{C}_1)\cup\mathcal{P}_b(\mathcal{C}_2)$ are distinct and relatively prime to $b-1$. To simplify our arguments, we define $P$ to be the product of the primes in $\mathcal{P}(\mathcal{C}_1)\cup\mathcal{P}(\mathcal{C}_2)$, and we let $r_{(0,0,0}=m_{(0,0,0)}=p_{(0,0,0)}=1$. 

Notice that if $2$ divides $b^m-1$, then $2$ divides $b-1$, and, hence, since all moduli of $\mathcal{C}_1\cup\mathcal{C}_2$ are greater than $1$, $2\not\in\mathcal{P}_b(\mathcal{C}_1)\cup\mathcal{P}_b(\mathcal{C}_2)$. If $3$ divides $P$, then $3$ is a primitive prime divisor of $b^m-1$ for some integer $m$. Thus, $b$ is not divisible by $3$ and by Fermat's little theorem, $b^2-1$ is divisible by $3$. It follows that $m\in\{1,2\}$, and since all moduli of $\mathcal{C}_1$ and $\mathcal{C}_2$ are greater than $1$, we deduce that $m=2$ and $b\equiv 2\pmod{3}$. Consequently, $3=p_{(1,0,0)}$ and $\gcd(3,b-1)=1$. 

Since $b^{2^q3^qq}\equiv b^{2^q3^q}\pmod{q}$, we have that if $q$ divides $P$, then $q=p_{(a,\overline{a},0)}$ for some $0\leq a\leq q$ and $0\leq\overline{a}\leq q$ with $(a,\overline{a},0)\neq(0,0,0)$. For the remainder of the proof, if $q$ divides $P$, then we define $a$ and $\overline{a}$ as mentioned here. Since all moduli of $\mathcal{C}_1\cup\mathcal{C}_2$ are greater than $1$, if $q$ divides $P$, then $q$ does not divide $b-1$, and therefore $\gcd(q,b-1)=1$. 

With the above observations, we can use the Chinese remainder theorem to let $\upsilon$ satisfy the congruences
    \begin{align*}
        \upsilon&\equiv-\alpha\pmod{2^q}\\
        \upsilon&\equiv(-\alpha-b^{(1,0,0)})\cdot b^{\alpha}\pmod{3^q},
    \end{align*}
let $\tau$ satisfy the congruences
    \begin{align*}
        \tau&\equiv -\alpha\pmod{2^q}\\
        \tau&\equiv
            \begin{cases}
                -\alpha\pmod{3^q}&\text{ if }3\nmid P\\ 
                (-\alpha-b^{-r_{(1,0,0)}})\cdot b^{\alpha}\pmod{3^q}&\text{ if }3\mid P
            \end{cases}\\
        \tau&\equiv
            \begin{cases}
            -\alpha\pmod{q}&\text{ if }q\nmid P\\
            (-\alpha-b^{-r_{(a,0,0)}})\cdot b^{\alpha}\pmod{q}&\text{ if }\overline{a}=0,3\nmid P\text{ and }q\mid P\\
            (-\alpha-b^{-r_{(a,0,0)}})\cdot b^{-\upsilon}\pmod{q}&\text{ if }\overline{a}=0,3\mid P\text{ and }q\mid P\\(-\alpha+b^{-r_{(a,\overline{a},0)}})\cdot b^{\alpha}\pmod{q}&\text{ if }\overline{a}\geq1,3\nmid P\text{ and }q\mid P\\
            (-\alpha+b^{-r_{(a,\overline{a},0)}})\cdot b^{-\upsilon}\pmod{q}&\text{ if }\overline{a}\geq1,3\mid P\text{ and }q\mid P, \end{cases}\\
    \end{align*}
let $\mathcal{T}$ satisfy the congruences
    \begin{equation*}
    \begin{aligned}\label{eq:brierfirstcongruences}
        \mathcal{T}&\equiv-\alpha\pmod{2^q}\\
        \mathcal{T}&\equiv-\alpha\pmod{b-1}\\
        \mathcal{T}&\equiv(-\alpha-b^{-r_{(j,0,j')}})\cdot b^{-\tau}\pmod{p_{(j,0,j')}}&\text{ for }0\leq j\leq q\text{ and }0\leq j'\leq 1 \\
        &&\text{ with }p_{(j,0,j')}\not\in\{3,q\}\\
        \mathcal{T}&\equiv(-\alpha+b^{-r_{(j,\overline{j},j')}})\cdot b^{-\tau}\pmod{p_{(j,\overline{j},j')}}&\text{ for }0\leq j\leq q, 1\leq\overline{j}\leq q\\
        &&\text{ and }0\leq j'\leq 1\\
        &&\text{ with }p_{(j,\overline{j},j')}\not\in\{3,q\},
    \end{aligned}
    \end{equation*}
and finally, we let $T$ satisfy the congruences 
    \begin{align*}
        T&\equiv\mathcal{T}\pmod{(b-1)\cdot P/\gcd(3q,P)}\\
        T&\equiv\tau\pmod{2^q3^qq}.
    \end{align*}
Then we can write $T=2^q3^qqh+\tau$ for some integer $h$, and since $b^{2^q3^qqh}\equiv 1\pmod{p_{(j,\overline{j},j')}}$ for each $0\leq j\leq q$, $0\leq\overline{j}\leq q$ and $0\leq j'\leq 1$,  when $p_{(j,0,j')}\not\in\{3,q\}$, we have that  
    \begin{align*}
        T\cdot b^T+\alpha
        &=T\cdot b^{(2^q3^qqh)}\cdot b^{\tau}+\alpha\\
        &\equiv(-\alpha-b^{-r_{(j,0,j')}})\cdot b^{-\tau}\cdot b^{(2^q3^qqh)}\cdot b^{\tau}+\alpha\pmod{p_{(j,0,j')}}\\
        &\equiv-b^{-r_{(j,0,j')}}\pmod{p_{(j,0,j')}}.       
    \end{align*}
Similarly, for $0\leq j\leq q$, $1\leq\overline{j}\leq q$ and $0\leq j'\leq 1$, when $p_{(j,\overline{j},j')}\not\in\{3,q\}$, we have that  
    \begin{align*}
        T\cdot b^T+\alpha
        &=T\cdot b^{(2^q3^qqh)}\cdot b^{\tau}+\alpha\\
        &\equiv(-\alpha+b^{-r_{(j,\overline{j},j')}})\cdot b^{-\tau}\cdot b^{(2^q3^qqh)}\cdot b^{\tau}+\alpha\pmod{p_{(j,\overline{j},j')}}\\
        &\equiv b^{-r_{(j,\overline{j},j')}}\pmod{p_{(j,\overline{j},j')}}. 
    \end{align*}

Suppose $3$ divides $P$. Then $3=p_{(1,0,0)}$ and $T\equiv\tau\pmod{2^q3^qq}$ imply $T\equiv -\alpha\pmod{2^q}$. Thus, $T=2^q\ell-\alpha$ for some integer $\ell$, and since $b^{2^q\ell}\equiv 1\pmod{p_{(1,0,0)}}$, we have that
    \begin{align*}
        T\cdot b^T+\alpha
        &=T\cdot b^{2^q\ell}b^{-\alpha}+\alpha\\
        &\equiv (-\alpha-b^{-r_{(1,0,0)}})\cdot b^{\alpha}\cdot b^{2^q\ell}b^{-\alpha}+\alpha\pmod{p_{(1,0,0)}}\\
        &\equiv-b^{-r_{(1,0,0)}}\pmod{p_{(1,0,0)}}.
    \end{align*}

Next, suppose $3$ divides $P$ and $q$ divides $P$. Then $q=p_{(a,\overline{a},0)}$ and $T\equiv\tau\pmod{2^q3^qq}$ imply $T\equiv\upsilon\pmod{2^q3^q}$. Thus, $T=2^q3^qh+\upsilon$ for some integer $h$. Since $b^{2^q3^qh}\equiv 1\pmod{p_{(a,\overline{a},0)}}$, if $\overline{a}=0$, we have that
    \begin{align*}
        T\cdot b^T+\alpha
        &= T\cdot b^{2^q3^qh}b^{\upsilon}+\alpha\\
        &\equiv (-\alpha-b^{-r_{(a,0,0)}})\cdot b^{-\upsilon}\cdot b^{2^q3^qh}b^{\upsilon}+\alpha\pmod{p_{(a,0,0)}}\\
        &\equiv-b^{-r_{(a,0,0)}}\pmod{p_{(a,0,0)}},
    \end{align*}
and if $\overline{a}\geq 1$, then
    \begin{align*}
        T\cdot b^T+\alpha
        &= T\cdot b^{2^q3^qh}b^{\upsilon}+\alpha\\
        &\equiv (-\alpha+b^{-r_{(a,\overline{a},0)}})\cdot b^{-\upsilon}\cdot b^{2^q3^qh}b^{\upsilon}+\alpha\pmod{p_{(a,0,0)}}\\
        &\equiv b^{-r_{(a,\overline{a},0)}}\pmod{p_{(a,0,0)}},
    \end{align*}

Lastly, suppose $3$ does not divide $P$ and $q$ divides $P$. Then $p_{(a,\overline{a},0)}=q$, and $T\equiv\tau\pmod{2^q3^qq}$ implies $T\equiv -\alpha\pmod{2^q3^q}$. Thus, $T=2^q3^qk-\alpha$ for some integer $k$. Since $b^{2^q3^qk}\equiv 1\pmod{p_{(a,\overline{a},0)}}$, if $\overline{a}=0$, we have that
\begin{align*}
    T\cdot b^T+\alpha
    &\equiv T\cdot b^{2^q3^qk}b^{-\alpha}+\alpha\pmod{p_{(a,0,0)}}\\
    &\equiv (-\alpha-b^{-r_{(a,0,0)}})\cdot b^{\alpha}\cdot b^{2^q3^qk}b^{-\alpha}+\alpha\pmod{p_{(a,0,0)}}\\
    &\equiv-b^{-r_{(a,0,0)}}\pmod{p_{(a,0,0)}}
\end{align*}
and if $\overline{a}\geq 1$, then 
    \begin{align*}
        T\cdot b^T+\alpha
        &\equiv T\cdot b^{2^q3^qk}b^{-\alpha}+\alpha\pmod{p_{(a,\overline{a},0)}}\\
        &\equiv (-\alpha+b^{-r_{(a,\overline{a},0)}}\cdot b^{\alpha})\cdot b^{2^q3^qk}b^{-\alpha}+\alpha\pmod{p_{(a,\overline{a},0)}}\\
        &\equiv b^{-r_{(a,\overline{a},0)}}\pmod{p_{(a,\overline{a},0)}}.
    \end{align*}

In all cases, if $n$ satisfies a congruence in $\mathcal{C}_1$, then $(Tb^T+\alpha)\cdot b^n+1$ is divisible by some prime in $\mathcal{P}(\mathcal{C}_1)$ and if $n$ satisfies a congruence in $\mathcal{C}_2$, then $(Tb^T+\alpha)\cdot b^n-1$ is divisible by some prime in $\mathcal{P}(\mathcal{C}_2)$. The condition that $T\equiv\mathcal{T}\pmod{b-1}\equiv -\alpha\pmod{b-1}$ ensures that $(Tb^T+\alpha)+1\equiv 1\pmod{b-1}$ and $(Tb^T+\alpha)-1\equiv -1\pmod{b-1}$, and, therefore, $\gcd((Tb^T+\alpha)+1,b-1)=\gcd((Tb^T+\alpha)-1,b-1)=1$. Since $\mathcal{C}_1$ and $\mathcal{C}_2$ are each covering systems, this finishes the proof.
\end{proof}

Our final corollary generalizes the results of Berrizbeitia et al. for integers $b\geq 2$ with $b$ not a power of $2$. When $b>2$, the corollary follows immediately from Theorem~\ref{thm:brier} and the covering systems provided in \eqref{eq:firstcovering} and \eqref{eq:lastcovering} by letting $\alpha\in\{-1,1\}$.

\begin{corollary}
For any integer $b\geq 2$ with $b+1$ not a power of $2$, there are infinitely many $b$-Brier numbers in the sequence of $b$-Cullen numbers and the sequence of $b$-Woodall numbers.
\end{corollary}

We end this section by providing a covering system satisfying the conditions on $\mathcal{C}$ in Theorem~\ref{thm:notbrierwithout2} and $\mathcal{C}_2$ in Theorem~\ref{thm:brier}. Using notation of Harrington \cite{harrington}, we let $([A_1,\ldots,A_v],[B_1,\ldots,B_v])$ represent the congruence $A\pmod{B}$, where $B=\text{lcm}(B_1,\ldots,B_v)$ and $A$ satisfies the congruences $A\equiv A_j\pmod{B_j}$ for each $1\leq j\leq v$. As an example, the covering system given in \eqref{eq:firstcovering} can be written as 
\[\{([2^{j-1}-1],[2^j]):1\leq j\leq q-1\}\cup\{([j,2^m-1],[q,2^j]):0\leq j\leq q-1\}.\]
Using this notation, the following provides a covering system satisfying conditions on $\mathcal{C}$ in Theorem~\ref{thm:notbrierwithout2} and $\mathcal{C}_2$ in Theorem~\ref{thm:brier}:
    \begin{align}\label{eq:lastcovering}
    \begin{split}
        &\{([3^{\overline{j}-1}],[3^{\overline{j}}]):1\leq\overline{j}\leq q\}\\
        \cup&\{([2^{j-1},2\cdot3^{\overline{j}-1}],[2^j,3^{\overline{j}}]):1\leq j\leq q, 1\leq\overline{j}\leq q\}\\
        \cup&\{([0,2\cdot3^{\overline{j}-1},j],[2^j,3^{\overline{j}},q]):1\leq j\leq q, 1\leq\overline{j}\leq q\}\\
        \cup&\{([0,\overline{j}],[3^{\overline{j}},q]):1\leq\overline{j}\leq q\}
    \end{split}
    \end{align}

\end{document}